\documentclass{amsart}

\usepackage{verbatim}
\usepackage{amsmath,amsthm,amssymb}
\usepackage[mathscr]{eucal}
\usepackage{enumerate}

\usepackage{eufrak}
\usepackage{graphicx}
\input{epsf}
\newtheorem{theorem}{Theorem}
\newtheorem{corollary}[theorem]{Corollary}
\newtheorem{proposition}[theorem]{Proposition}

\newtheorem{lemma}[theorem]{Lemma}

\theoremstyle{definition}
\newtheorem{definition}[theorem]{Definition}
\newtheorem{remark}[theorem]{Remark}
\newtheorem{example}[theorem]{Example}

\newcommand{\POVM}{positive operator-valued measure }
\newcommand{\wotlim}{\ensuremath{\hbox{\sc wot}\,\lim}}
\newcommand{\wotclos}{\ensuremath{\hbox{\sc wot-}\text{cl }}}

\begin{document}

\title[$C^*$-Extreme Maps]{On $C^*$-Extreme Maps and $*$-Homomorphisms of a Commutative
$C^*$-Algebra\footnote{T\lowercase{his paper constitutes a part of
the author's }P\lowercase{h}D. \lowercase{thesis at the
}U\lowercase{niversity of
}N\lowercase{ebraska-}L\lowercase{incoln}}}

\author[M. C. Gregg]{Martha Case Gregg}
\address{Dept. of Mathematics\\
Augustana College\\2001 South Summit Avenue\\ Sioux Falls, SD 57197}
\email{martha.gregg@augie.edu}

\begin{abstract}
The generalized state space of a commutative $C^*$-algebra, denoted
$S_\mathcal{H}(C(X))$, is the set of positive unital maps from
$C(X)$ to the algebra $\mathcal{B(H)}$ of bounded linear operators
on a Hilbert space $\mathcal{H}$.  $C^*$-convexity is one of several
non-commutative analogs of convexity which have been discussed in
this context. In this paper we show that a $C^*$-extreme point of
$S_\mathcal{H}(C(X))$ satisfies a certain spectral condition on the
operators in the range of the associated positive operator-valued
measure.  This result enables us to show that $C^*$-extreme maps
from $C(X)$ into $\mathcal{K}^+$, the algebra generated by the
compact and scalar operators, are multiplicative.  This generalizes
a result of D. Farenick and P. Morenz.  We then determine the
structure of these maps.

\end{abstract}

 \maketitle

Several non-commutative analogs of convexity have appeared in the
literature including CP-convexity \cite{FujimotoGeNaThfoC*Al98} and
matrix convexity \cite{EffrosWinlkerMaCoOpAnBiHBTh}, as well as
$C^*$-convexity \cite{FarenickMorenzC*ExPoGeStSpC*Al},
\cite{HopenwasserMoorePaulsenC*ExPo}, which is the topic of this
paper.  In \cite{HopenwasserMoorePaulsenC*ExPo}, Hopenwasser, Moore,
and Paulsen characterized operators which are $C^*$-extreme in the
unit ball of $\mathcal{B(H)}$ and obtained results about other
$C^*$-convex sets and their extreme points.  In
\cite{FarenickMorenzC*ExPoGeStSpC*Al}, Farenick and Morenz extend
the idea of $C^*$-convexity to the space of completely positive maps
on a $C^*$-algebra.  They show that $C^*$-extreme maps with their
range in $\mathcal{K}^+$ are also extreme (in the classical sense)
and obtain a characterization of $C^*$-extreme maps on a commutative
$C^*$-algebra which have their range in $M_n$, the $C^*$-algebra of
$n \times n$ complex matrices.  Subsequently, Zhou
\cite{ZhouC*ExPoSpCoPoMa} gave two necessary and sufficient
conditions for a completely positive map to be $C^*$-extreme,  and
described the structure of $C^*$-extreme maps with range in $M_n$.
 The main results presented here are Theorem~\ref{nec}, which gives a
necessary condition for a map $\phi:C(X) \rightarrow \mathcal
{B(H)}$ to be $C^*$-extreme, and Theorem~\ref{cmpct}, which then
shows that a positive unital map $\phi:C(X)\rightarrow\mathcal{K}^+$
is $C^*$-extreme if and only if it is multiplicative.  We then
determine the structure of such maps.

The author wishes to express her gratitude to Professor David Pitts
for many lively conversations which led to significant improvements
in the paper.

Throughout, let $X$ be a compact Hausdorff space, $C(X)$ the
$C^*$-algebra of continuous functions on $X$, and $\mathcal{H}$ a
Hilbert space.
\begin{definition}
The \emph{generalized state space} of $C(X)$ is
    $$S_\mathcal{H}(C(X)) = \{\phi:C(X) \rightarrow \mathcal{B(H)}
    \,|\, \phi \text{ is positive and } \phi(1) = I\}.$$
\end{definition}
Note that in the case of a non-commutative $C^*$-algebra
$\mathfrak{A}$, the generalized state space $S_\mathcal{H}(A)$ is
the set of \emph{completely positive} unital maps.  However, for a
commutative $C^*$-algebra, every positive map is also completely
positive \cite[Theorem 4]{StinespringPoFuC*Al}.  If $\mathcal{H}=
\mathbb{C}$, the generalized state space
$S_\mathbb{C}(\mathfrak{A})$ coincides with the classical state
space of $\mathfrak{A}$.
\begin{definition}
We say that  $\phi,\psi \in S_\mathcal{H}(C(X))$ are \emph{unitarily
equivalent}, and write $\phi \sim \psi$, if there is a unitary $u
\in \mathcal{B(H)}$ such that $\phi(f) = u^*\psi(f)u$ for every $ f
\in C(X)$.
\end{definition}
\begin{definition}If $\phi, \psi_1, \ldots
\psi_n \in S_\mathcal{H}(C(X))$ and $t_1,\ldots t_n \in
\mathcal{B(H)}$ are invertible with $t_1^*t_1 + \ldots + t_n^*t_n =
I$, then we say
    $$\phi(f) = t_1^* \psi_1(f)t_1 + \ldots + t_n^* \psi_n(f)t_n \text{ for every } f\in C(X), $$
is a \emph{proper $C^*$-convex combination}.  We call a map $\phi
\in S_\mathcal{H}(C(X))$ \emph{$C^*$-extreme} if, whenever $\phi$ is
written as a proper $C^*$-convex combination of $\psi_1,\,
\ldots,\,\psi_n$, then $\psi_j \thicksim \phi$ for each $j = 1,
\ldots,n$.
\end{definition}


We begin with a discussion of $\mathcal{B(H)}$-valued measures,
which closely follows the development given in Paulsen
\cite{PaulsenCoBoMaOpAl}. These operator valued measures play a key
role in the proof of Theorem~\ref{nec}, below.  Given a bounded
linear map $\phi:C(X) \rightarrow \mathcal{B(H)}$ and vectors $x,y
\in \mathcal {H}$, the bounded linear functional
    $$f \mapsto \langle \phi(f) x, y \rangle$$
corresponds to a unique regular Borel measure $\mu_{x,y}$ on $X$
such that
    $$\int_X f d\mu_{x,y} := \langle \phi(f) x,y \rangle \text{ for any }f \in C(X).$$
Denote the $\sigma$-algebra of Borel sets of $X$ by $\mathcal{S}$.
For a set $B \in \mathcal{S}$, the sesquilinear form
    $$(x,y) \mapsto \mu_{x,y}(B)$$
then determines an operator $\mu(B)$.  Thus we obtain an
\emph{operator-valued measure} $\mu:\mathcal{S}
\longrightarrow\mathcal{B(H)}$ which is:
    \begin{enumerate}
        \item \emph{weakly countably additive},  i.e., if $\{B_i\}_{i=1}^\infty \subseteq \mathcal{S}$ are
        pairwise disjoint, and $B = \bigcup_{i=1}^\infty B_i$ then
        $$\langle \mu(B)x,y \rangle = \sum_{i=1}^\infty \langle
        \mu(B_i) x, y\rangle \text{ for every } x,y \in
        \mathcal{H}.$$
        \item \emph{bounded}, i.e., $\|\mu \|:= \sup \{\|\mu(B) \|:B \in \mathcal{S}\} <
        \infty$.
        \item \emph{regular}, i.e., for each pair of vectors $x$ and
        $y$ in $\mathcal{H}$, the complex measure $\mu_{x,y}$ is
        regular.
    \end{enumerate}
Furthermore, this process works in reverse: given a regular bounded
operator-valued measure $\mu:\mathcal{S} \longrightarrow
\mathcal{B(H)}$, define Borel measures
    $$\mu_{x,y}(B):= \langle \mu(B)x,y \rangle$$
for each $x,y \in \mathcal{H}$.  Then the operator $\phi(f)$ is
uniquely defined by the equations
    $$\langle \phi(f)x,y \rangle := \int_X f d\mu_{x,y};$$
the map $\phi:C(X) \longrightarrow \mathcal{B(H)}$ is then seen to
be bounded and linear.  This construction shows that each operator
valued-measure gives rise to a unique bounded linear map, and
vice-versa.  The following proposition summarizes properties shared
by operator valued-measures and their associated linear maps.  We
will be most concerned with parts (2) and (4) of Proposition~\ref{paulsen4.5}; part (4) is, of
course, the Spectral Theorem.

\begin{proposition}\label{paulsen4.5}
\cite[Proposition 4.5]{PaulsenCoBoMaOpAl} Given an operator valued
measure $\mu$ and its associated linear map $\phi$,
\begin{enumerate}
    \item $\phi$ is self-adjoint if and only if $\mu$ is self-adjoint,
    \item $\phi$ is positive if and only if $\mu$ is positive,
    \item $\phi$ is a homomorphism if and only if $\mu(B_1 \cap B_2)
    = \mu(B_1)\mu(B_2)$ for all $B_1, B_2~\in~\mathcal{S}$,
    \item $\phi$ is a $*$-homomorphism if and only if $\mu$ is spectral
    (i.e., projection-valued).
\end{enumerate}
\end{proposition}

We note the following important features of positive operator valued
measures, and their relationship to the associated positive maps.

    \textbf{(1)} Let $\mathfrak{F}(X) = \{f:X \rightarrow \mathbb{C}\,|\,\, f \text{
    is a bounded Borel measurable function}\}$. If $\phi:C(X)
    \rightarrow \mathcal{B(H)}$ is a positive map, we may use the
    corresponding \POVM  to extend $\phi$ to a map
    $\tilde{\phi}:\mathfrak{F}(X) \rightarrow \mathcal{B(H)}$ by
    defining
        $$\tilde{\phi}(f) = \int_X f d\mu_\phi,$$
    for every $f \in \mathfrak{F}(X) $.  The measure $\mu_\phi$ may then
    be viewed as the restriction of $\tilde{\phi}$ to the characteristic
    functions of Borel sets.  For
    simplicity, we will simply write $\phi$, rather than $\tilde{\phi}$,
    and use the notations $\mu_\phi(F)$ and $\phi(\chi_{_F})$
    interchangeably.

    \textbf{(2)} A positive unital map $\phi \in S_H(C(X))$ is $C^*$-extreme if and only if
    the associated operator-valued measure $\mu_\phi$ is $C^*$-extreme.  (Here,
    a positive operator-valued measure $\mu_\phi$ is called $C^*$-extreme if,
    whenever $\mu_\phi$ is written
        $$\mu_\phi = t_1^* \mu_1 t_1 + ... + t_n^* \mu_n t_n,$$
    where $\sum_{j = 1}^n t_j^* t_j = I$ and each $\mu_j$ is a
    positive operator-valued measure, then $\mu_j \sim\mu_\phi$ for
    each $j = 1,...n$.)

    \textbf{(3)} Finally, if $\phi:C(X) \rightarrow \mathcal{B(H)}$ is a positive bounded linear
        map, and $\mu_\phi$ the associated operator-valued measure, then for each Borel set $F~\subseteq~X$,
         $\mu_\phi(F)~\in~\wotclos\phi(C(X))$, the weak operator topology closure of
         $\phi(C(X))$.  The proof of this fact requires some care,
         because while $\phi(C(X))$ is an operator space, it is
         not generally an algebra.

    \begin{proof}
        Let $G \subseteq X$ be an open set.  Then a basic WOT-open set in
        $\mathcal{B(H)}$ centered at $\phi(\chi_{_G})$ has the form:
            $$\mathcal{O} = \{T \in \mathcal{B(H)}\,:\,|\langle (T - \phi(\chi_{_G}))x_i,y_i \rangle | <
            \varepsilon
            \text{ for }i = 1 \ldots n\},
            $$
    where $x_i, y_i \in \mathcal{H}$ and $\varepsilon > 0$.  We wish
    to show that for any such open set there is a function $f \in C(X)$
    with $\phi(f) \in \mathcal{O}$.  For each $j$, we can write
        $$\mu_{x_j,y_j} = \mu_{j,1}- \mu_{j,2} + i (\mu_{j,3} - \mu_{j,4}),$$
    where each $\mu_{j,k}$ is a positive measure.  Since each of
    these measures is regular, we may choose compact sets $K_{j,k} \subseteq G$
    for $j = 1, \ldots, n$ and $k = 1, \ldots, 4$ such that
        $$\mu_{j,k}(G \setminus K_{j,k}) < \frac{\varepsilon}{4}.$$
    Then, setting
        $$\displaystyle{K = \bigcup_{j = 1}^n \bigcup_{k = 1}^4 K_{j,k}}$$
    we have, for $j = 1,\ldots, n$,
        $$|\mu_{x_j,y_j}(G \setminus K)| \leq |\mu_{j,1}(G \setminus K)|+ \cdots
        + |\mu_{j,4}(G \setminus K)| < \varepsilon.$$
    Urysohn's Lemma now guarantees the existence of a continuous
    function $f:X \rightarrow [0,1]$ with $f|_K = 1$ and $f|_{G^C} =
    0$.  Hence, for each $j = 1, \ldots, n$,
        \begin{align*}
        |\langle (\phi(f) - \mu_\phi(K))x_j, y_j \rangle | & =
        \biggl| \int_X (f - \chi_{_K}) d\mu_{j,1} - \int_X (f - \chi_{_K})
        d\mu_{j,2}\\
        & \hspace{.3 in} + i\left( \int_X (f - \chi_{_K})d\mu_{j,1}
            - \int_X (f - \chi_{_K}) d\mu_{j,4}\right)\biggr|\\
        & \leq \int_X \chi_{_{G\setminus K}} d\mu_{j,1} + \cdots +
        \int_X \chi_{_{G\setminus K}} d\mu_{j,4}\\
        & < \varepsilon.
        \end{align*}
    Therefore $\phi(f) \in
    \mathcal{O}$, as required; hence $\phi(f) \in
    \wotclos\phi(C(X))$.

    Now let
        $$\mathcal{F} = \left\{F \subseteq X\, : \text{ F is a Borel set and
        } \phi(\chi_{_F}) \subseteq \wotclos\phi(C(X)) \right\}.$$
    We will prove that $\mathcal{F}$ is a $\sigma$-algebra
    containing the Borel sets, and hence that $\mathcal{F} =
    \mathcal{S}$.  Our discussion above shows that $\mathcal{F}$ contains every open set of $X$.
    Suppose that $\{B_i\}$ is a countable family of sets in $\mathcal{F}$ and set $B~=~\bigcup_{i = 1}^\infty B_i$.
    Assume without loss of generality that $\{B_i\}$ are a
    disjoint family. Then, since
    $\mu_\phi$ is weakly countably additive,
        $$\langle \mu_\phi(B)x, y \rangle = \sum_{i = 1}^\infty
        \langle \mu_\phi(B_i)x, y \rangle$$
    for any $x, y \in \mathcal{H}$.  That is
        $$\displaystyle{\mu_\phi(B) = \wotlim_N \mu_\phi\biggl(\bigcup_{i=1}^N B_i\biggr)};  $$
    It follows that $B \in \mathcal{F}$.  Furthermore, if $F \in \mathcal{F}$, then
        $$\phi(\chi_{_{F^C}}) = \phi(1 -\chi_{_F}) = I - \phi(\chi_{_F}),$$
    so that $F^C \in \mathcal{F}$ also.
    Therefore $\mathcal{F}$ is the
    $\sigma$-algebra of Borel sets of $X$.
    \end{proof}
    Thus, if the range of $\phi$ is contained
    in a $C^*$-subalgebra $\mathcal{A}$ of $\mathcal{B(H)}$, then the
    range of $\mu_{\phi}$ is contained in the weak operator topology
    closure of $\mathcal{A}$, i.e, $\mathcal{A}^{\prime\prime}$.

We can now prove the following theorem, which gives a necessary
condition for a positive map $\phi$ on a commutative $C^*$-algebra
(or equivalently its associated  \POVM) to be $C^*$-extreme.

\begin{theorem} \label{nec}
Let $X$ be a compact Hausdorff space, and $\phi:C(X) \longrightarrow
\mathcal{B(H)}$ a unital, positive map. Denote by $\mu_{\phi}$ the
unique  \POVM  associated to $\phi$. If $\phi$ is $C^*$-extreme,
then for every Borel set $F \subset X$, either
    \begin{description}
        \item[(1)] $\mu_{\phi}(F)$ is a projection, in which case $\mu_\phi(F) \in \phi(C(X))'$, or
        \item[(2)] $\sigma(\mu_{\phi}(F))= [0,1]$.
    \end{description}
Moreover, if (2) occurs and $\mu_\phi(F)$ has an eigenvalue in
$(0,1)$, then the point spectrum of $\mu_{\phi}(F)$ must contain
$(0,1)$.
\end{theorem}

\begin{proof}
Suppose there is a Borel set $F \subseteq X$ so that $\mu_{\phi}(F)$
is not a projection and $\sigma(\mu_{\phi}(F))\neq [0,1]$.  We will
show that $\phi$ is not $C^*$-extreme by constructing a proper
$C^*$-convex combination
    $$t_1^*\psi_1t_1 + t_2^*\psi_2t_2 = \phi$$
in which $\psi_1$ and $\psi_2$ are not unitarily equivalent to
$\phi$. Choose $x~\in~(0,1)\setminus\sigma(\mu_{\phi}(F))$ and
let $(a,b)$ be the largest open subinterval of $(0,1)$ which
contains $x$ but does not intersect $\sigma(\mu_{\phi}(F))$.  To be
precise, let
    $$(a,b) = \bigcup\{(\alpha, \beta)\subseteq (0,1): x \in
    (\alpha, \beta), \, (\alpha,\beta) \cap \sigma(\mu_\phi(F)) =
    \varnothing \}$$
Note that this choice of the interval $(a,b)$ insures that at least
one of the pair $\{a,b\}$ is in $\sigma(\mu_{\phi}(F))$.  In
particular, if $a >0$ then $a \in \sigma(\mu_{\phi}(F))$ and if $b <
1$ then $b \in \sigma(\mu_{\phi}(F))$. Choose $s_1 \in
\left(\frac{1}{4},\frac{1}{2}\right)$ with $s_1 >\frac{1}{2} \left(\frac{a - ab}{b - ab}\right)$,
and set $s_2 = 1 - s_1$. For $k = 1, 2$, define
    $$Q_k =\frac{1}{2}\mu_{\phi}(F) + s_k \mu_{\phi}(F^C) = s_kI +
    (\frac{1}{2}- s_k)\mu_{\phi}(F).$$
Note that $0 \not\in \sigma(Q_k) = s_k + (\frac{1}{2}- s_k)
\sigma(\mu_{\phi}(F))$, so that both $Q_k$'s are invertible. Now
define new  positive operator-valued measures $\mu_1$ and $\mu_2$ by
    $$\mu_k(B) = Q_k^{- \frac{1}{2}} \left(\frac{1}{2}\mu_{\phi}(B \cap
    F) + s_k\mu_{\phi}(B \cap F^C)\right)Q_k^{- \frac{1}{2}},$$
where $B$ is any Borel set of $X$. Observe that each of the
$\mu_k$'s is a  \POVM  with $\mu_k(X) = I$. Next, define $t_k =
Q_k^{\frac{1}{2}}$, for $k = 1,2$. Then, for any Borel set $B$ of
$X$,

\begin{align*}
    t_1^* \mu_1(B)t_1 + t_2^* \mu_2(B) t_2 &= \frac{1}{2}\mu_{\phi}(B \cap F) + s_1\mu_{\phi}(B
    \cap F^C)\\
    & \hspace{.7 in} + \frac{1}{2}\mu_{\phi}(B \cap F) + s_2 \mu_{\phi}(B \cap F^C)\\
    & = \mu_{\phi}(B) .
\end{align*}

Each $t_k$ is invertible and
    $$t_1^*t_1 + t_2^*t_2 = Q_1 + Q_2 = \mu_{\phi}(F) + \mu_{\phi}(F^C)=
    I.$$
Thus $t_1^* \mu_1t_1 + t_2^* \mu_2 t_2$ is a proper $C^*$-convex
combination of $\mu_1$ and $\mu_2$.

It is still necessary to show that $\mu_{\phi}$ is not unitarily
equivalent to at least one of $\mu_1$ or $\mu_2$. For $k = 1,2$, set $g_k(t) = [s_k +
(s_k - \frac{1}{2})t]^{-\frac{1}{2}}$.   As each $g_k$ is continuous on
$[0,1]$, and $Q_k^{-\frac{1}{2}}= g_k(\mu_{\phi}(F))$,
$Q_k^{-\frac{1}{2}}$ commutes with $\mu_{\phi}(F)$. Thus, for $k =
1, 2$, we have

\begin{align*}
    \mu_k(F) &= Q_k^{-\frac{1}{2}}
    \left(\frac{1}{2}\mu_{\phi}(F)\right) Q_k^{-\frac{1}{2}}\\
    &=\frac{1}{2}\mu_{\phi}(F)\left(s_k I + \left(\frac{1}{2}-
    s_k\right)\mu_{\phi}(F)\right)^{-1}.
\end{align*}

Let $f_k(t) = \frac{1}{2}t\left(s_k +
\left(\frac{1}{2}-s_k\right)t\right)^{-1} $.  Observe that each
$f_k$ is continuous on $[0,1]$, and that $\mu_k(F) =
f_k(\mu_{\phi}(F))$.  Therefore, by the spectral mapping theorem,
$\sigma(\mu_k(F))=f_k(\sigma(\mu_{\phi}(F)))$.  It is easy to check
that for $t \in (0,1)$, $t < f_1(t) < 1$, while $0 < f_2(t) < t$,
and that both $f_k$'s are strictly increasing. In addition, since
$s_1
> \frac{1}{2} \left(\frac{a - ab}{b - ab}\right)$, if $a>0$,
    $$a < f_1(a) = \frac{1}{2}a\left(\frac{1}{s_1 +
    (\frac{1}{2}-s_1)a}\right)<\frac{a}{\frac{a(1-b)}{b(1 -a)}(1-a)
    +a}=b \leq f_1(b).$$
Consider the following two cases:

\begin{description}
    \item[Case (i)] $a \neq 0$.  In this case $a \in
    \sigma(\mu_{\phi}(F))$.
  Thus $f_1(a) \in \sigma(\mu_1(F))$, but since $f_1(a)\in(a,b)$,
  $f_1(a) \notin
    \sigma(\mu_{\phi}(F))$. This shows that $\sigma(\mu_{\phi}(F)) \neq
    \sigma(\mu_1(F))$; therefore $\mu_{\phi}$ and $\mu_1$ are not
    unitarily equivalent.

    \item[Case (ii)] $a = 0$.  In this case, $b<1$ and $b \in
    \sigma(\mu_{\phi}(F))$. As $a = 0< f_2(b) < b$, we have $f_2(b) \in
    \sigma(\mu_2(F)) \setminus\sigma(\mu_{\phi}(F))$.  In this case, $\mu_2$ is
     not unitarily equivalent to $\mu_{\phi}$.
\end{description}

Let $\psi_k$ be the positive map determined by $\mu_k$.  Then $\phi
= t_1^* \psi_1 t_1 + t_2^* \psi_2 t_2$; this is a proper
$C^*$-convex combination of $\psi_1$ and $\psi_2$, where $\phi$ is
not unitarily equivalent to at least one of the maps $\psi_k$. Therefore, $\phi$ is
not $C^*$-extreme.

Now suppose that $\sigma(\mu_\phi(F)) = [0,1]$ and that $\sigma_{pt}(\mu_\phi(F))$
intersects $(0,1)$, but does not contain $(0,1)$.  It is not difficult to convince oneself that it is possible to choose $a,b \in (0,1)$ satisfying both
    \begin{enumerate}[(i)]
    \item $a < b < \frac{2a}{a+1}$, and
    \item exactly one of the pair $\{a,b\}$ is an eigenvalue.
    \end{enumerate}
Set $s_1 = \frac{1}{2}\left(\frac{a -ab}{b - ab}\right)$ and define positive operator-valued
measures $\mu_1$ and $\mu_2$ as above.  As in the previous computation, $\mu_1(F) = f_1(\mu_\phi(F))$.
As a result of our choice of $s_1$, $f_1(a) = b$.  Application of the Spectral Mapping Theorem then shows that either
    \begin{align*}
    b &\in \sigma_{pt}(\mu_1(F)) \setminus \sigma_{pt}(\mu_\phi(F)), \,\text{ or}\\
    b &\in \sigma_{pt}(\mu_\phi(F)) \setminus \sigma_{pt}(\mu_1(F)).
    \end{align*}
Since the point spectrum is
also a unitary invariant, and $\mu_\phi = t_1^*\mu_1t_1 + t_2^*\mu_2t_2$, this shows that $\phi$ is not
$C^*$-extreme.

Finally, we wish to show that any projection in the range of
$\mu_\phi$ must commute with $\phi(C(X))$.  Suppose that
$\mu_\phi(F)$ is a projection and choose $f \in C(X)$ with $0 \leq f
\leq 1$. Write
    $$f = \chi_{_F}f + (1 - \chi_{_F})f.$$
Then $\phi(\chi_{_F}f) \leq \mu_\phi(F)$, so these operators
commute. Similarly,
    $$\phi((1 - \chi_{_F})f) \leq \mu_\phi(X
    \setminus F) = I - \mu_\phi(F),$$
so that $\phi((1 - \chi_{_F})f)$ also commutes with $\mu_\phi(F)$.
Therefore $\phi(f)$ commutes with $\mu_\phi(F)$.  If $f$ is an
arbitrary continuous function, we can express $f$ as a linear
combination of positive functions with ranges in $[0,1]$.  Thus $f$
will commute with $\mu_\phi(F)$.
\end{proof}

In their paper of 1997 \cite{FarenickMorenzC*ExPoGeStSpC*Al},
Farenick and Morenz show that a positive map from a commutative
$C^*$-algebra into a matrix algebra $M_n$ is $C^*$-extreme if and
only if it is a $*$-homomorphism.  In view of the spectral condition
given by Theorem~\ref{nec}, a shorter proof is possible.

\begin{corollary}
\cite[Proposition 2.2]{FarenickMorenzC*ExPoGeStSpC*Al} Let $X$ be a
compact Hausdorff space and $\phi:C(X) \longrightarrow M_n$ a
positive map. Then $\phi$ is $C^*$-extreme if and only if it is a
$*$-homomorphism.
\end{corollary}

\begin{proof}
 It is already known that if $\phi$ is a representation
(i.e.,$*$-homomorphism), then $\phi$ is $C^*$-extreme
\cite[Proposition 1.2]{FarenickMorenzC*ExPoGeStSpC*Al}.  On the
other hand, if $\phi$ is not a representation, then the associated
\POVM $\mu_{\phi}$ is not a spectral measure.  In this case, there
is a Borel set $F \subset X$ for which $\mu_{\phi}(F)$ is not a
projection.  As $\mu_{\phi}(F)$ is an $n \times n$ matrix,
$\sigma(\mu_{\phi}(F))$ consists of at most $n$ isolated points. We
may therefore apply the theorem to conclude that $\phi$ is not
$C^*$-extreme.
\end{proof}

Note that in the proof of Theorem~\ref{nec}, $Q_k$,
$Q_k^{-\frac{1}{2}}$ and $t_k = Q_k^{\frac{1}{2}}$ are elements of
the $C^*$-algebra generated by $\mu_{\phi}(F)$.  As noted in the
remark preceding Theorem~\ref{nec}, the range of $\mu_{\phi}$ is
contained in the WOT-closure of the range of $\phi$. Thus we have
the following corollary:

\begin{corollary}
Let $\mathcal{M}\subseteq \mathcal{B(H)}$ be a von Neumann algebra,
$\phi:C(X) \longrightarrow \mathcal{M}$ a unital positive map, and
$\mu_{\phi}$ the positive operator-valued measure associated to
$\phi$. If $\phi$ fails to meet the spectral condition described in
Theorem~\ref{nec}, then $\phi$ can be written as a proper
$C^*$-convex combination
$$\phi = t_1^*\psi_1t_1+t_2^*\psi_2t_2,$$ where each $t_k \in
\mathcal{M}$, each $\psi_k:C(X) \longrightarrow \mathcal{M}$, and, for at least one choice of $k$,
$\psi_k$ is not unitarily equivalent to $\phi$ in
$\mathcal{B(H)}$.
\end{corollary}

We now consider an example of a $C^*$-extreme map which is not
multiplicative. The positive map $\phi$ defined below was considered
by Arveson \cite[p. 164]{ArvesonSuAlC*Al} as an example of an
extreme point in the generalized state space. Farenick and Morenz
\cite[Example 2]{FarenickMorenzC*ExPoGeStSpC*Al} subsequently showed
that $\phi$ is also a $C^*$-extreme point, although not a
homomorphism.  Consider the Hilbert spaces $L^2(\mathbb{T}, m)$,
where $m$ is normalized Lebesgue measure on $\mathbb{T}$, and $H^2$,
the classical Hardy space.  Let $P$ be the projection of
$L^2(\mathbb{T}, m)$ onto $H^2$. For a function $f \in
L^2(\mathbb{T}, m)$ denote by $M_f$ multiplication by $f$ and by
$T_f = PM_fP$ the Toeplitz operator for $f$.

\begin{example}\label{AFM}\cite{ArvesonSuAlC*Al},
\cite{FarenickMorenzC*ExPoGeStSpC*Al} Consider the representation
$\pi:C(\mathbb{T}) \longrightarrow \mathcal{B}(L^2(\mathbb{T},m))$
given by $\pi(f) = M_f$.  The spectral measure associated to $\pi$
is given by $\mu_{\pi}(B) = M_{\chi_{_B}}$, where $B \subseteq X$ is
a Borel set. Define a unital positive map
    $$\phi:C(\mathbb{T}) \longrightarrow \mathcal{B}(H^2)$$
by
    $$\phi(f) = PM_fP.$$
Since $\mu_{\pi}(B) = M_{\chi_{_B}}$, we have $\mu_{\phi}(B) =
PM_{\chi_{_B}}P = T_{\chi_{_B}},$ a Toeplitz operator. Thus
$\sigma(\mu_{\phi}(B)) = \sigma(T_{\chi_{_B}})$.  Since $\chi_{_B}$
is a real-valued $L^\infty$ function, $\sigma(T_{\chi_{_B}})$ is the
closed convex hull of the essential range of $\chi_{_B}$ \cite[p.
868]{HartmanWintnerSpToMa}. Therefore, if $\mu_{\phi}(B) \not \in
\{0,I\}$, then $\sigma(\mu_{\phi}(B)) = [0,1]$ . Thus, for any Borel
set $B \subseteq X$, either $\mu_{\phi}(B) = [0,1]$ or
$\mu_{\phi}(B)$ is a trivial projection; that is, $\phi$ satisfies
the conditions of the theorem.
\end{example}

Now let us consider the case of a unital positive map $\phi$ on a
commutative $C^*$-algebra $C(X)$ whose range is in $\mathcal{K}^+$,
the $C^*$-algebra generated by the compact operators and the
identity operator. In \cite[Proposition
1.1]{FarenickMorenzC*ExPoGeStSpC*Al} Farenick and Morenz show that
if such a map $\phi$ is $C^*$-extreme, then $\phi$ is also extreme.
It is possible, however, to say more. Theorem~\ref{nec} requires the
operators in the range of the positive operator-valued measure
$\mu_{\phi}$ either to be projections, or to have spectrum equal to
$[0,1]$. In contrast, the spectrum of a positive operator $K +
\alpha I \in \mathcal{K}^+$ must be a sequence of positive numbers
with a single limit point at $\alpha$.  This dichotomy suggests that
Theorem~\ref{nec} may give additional information about these maps.
In fact, both the result of Theorem~\ref{nec} (the spectral
condition on the operators in the range of $\mu_\phi$) and the
technique used in its proof, will be used below.  The result is
Theorem~\ref{cmpct}, which shows that such maps must be
multiplicative, and gives their structure.

In the succeeding lemma and theorem, let $q$ be the usual quotient
map $q:\mathcal{B(H)} \rightarrow \mathcal{B(H)/K(H)}$, and set
$\tau = q \circ \phi$. Then $\tau$ is a positive linear functional,
so there is a unique positive real-valued Borel measure $\mu_{\tau}$
on $X$ so that
    $$\tau(f) = \int_X f d\mu_{\tau} \text{ for every } f \in C(X).$$
For any function $f \in C(X)$, write
    $$\phi(f) = K_f + \tau(f)I,$$
where $K_f \in \mathcal{K}$ is a compact operator.

\begin{lemma}\label{taumult} Let $\phi:C(X) \rightarrow \mathcal{K}^+$ be unital, positive,
    and $C^*$-extreme.  Then the map $\tau $ is
    multiplicative.
\end{lemma}

\begin{proof}
    As in the proof of Theorem~\ref{nec}, we will prove the
    contrapositive.  Assume that $\tau$ is not multiplicative; then
    the support of $\mu_{\tau}$ must contain at least two distinct
    points, which we will call $s_1$ and $s_2$.  Let $N_1$ be a
    neighborhood of $s_1$ which does not contain $s_2$.  By
    Urysohn's Lemma, there exists a continuous function $f:X
    \rightarrow [0,1]$ such that $f(s_1) =1$ and $f|_{N_1^C} = 0$.

    Choose $\alpha$ and $\beta$ in $(0,1)$ with $\alpha > \beta$ and
    let
        \begin{align*}Q_1 &= \alpha \phi(f) + \beta \phi(1-f) = (\alpha -
        \beta)\phi(f) + \beta I, \text{ and}\\
        Q_2 &= ( 1 - \alpha) \phi(f) + (1 - \beta) \phi(1-f) = (\beta -
        \alpha)\phi(f) + (1 - \beta) I.
        \end{align*}
    Note that since $0 \leq f \leq 1$, the spectrum of $\phi(f)$ is
    contained in the closed unit interval.  Thus,
        \begin{align*}\sigma(Q_1) &\subseteq [\beta, \alpha],\text{
        and}\\
        \sigma(Q_2) &\subseteq [1 - \alpha, 1-\beta].
        \end{align*}
    So both $Q_j$'s are invertible positive operators.  Define maps
    $\psi_1$ and $\psi_2$ by
        \begin{align*}\psi_1(g) &= Q_1^{-\frac{1}{2}}\left[\alpha \phi(fg) +
        \beta \phi((1-f)g)\right] Q_1^{-\frac{1}{2}}, \text{ and}\\
        \psi_2(g) &= Q_2^{-\frac{1}{2}}\left[(1-\alpha) \phi(fg) +
        (1 -\beta) \phi((1-f)g)\right] Q_2^{-\frac{1}{2}}.
        \end{align*}
    Both $\psi_j$'s are positive, unital maps with ranges in $\mathcal{K}^+$.  Setting $t_j =
    Q_j^{\frac{1}{2}}$, we have
        \begin{align*}t_1^* \psi_1(g)t_1 + t_2^* \psi_2(g)t_2 &= \alpha \phi(fg)
        + \beta \phi(g - fg) + (1-\alpha)\phi(fg) + (1 -
        \beta)\phi(g - fg)\\
        & = \phi(fg) + \phi(g - fg)\\
        & = \phi(g), \,\, \text{ for every } g \in C(X).
        \end{align*}
    Since $t_1^*t_1 + t_2^*t_2 = I$, the above expression gives
    $\phi$ as a proper $C^*$-convex combination of $\psi_1$ and $\psi_2$.

    We now wish to show that $\psi_1$ and $\psi_2$ are not unitarily
    equivalent.  To this end, let $N_2$ be a neighborhood of $s_2$
    with $N_1 \cap N_2 = \emptyset$.  Then we may
    choose a continuous function $h:X \rightarrow [0,1]$
    with $h|_{N_2^C} = 0$ (i.e., supp $h \subseteq N_2$) and $h(s_2) = 1$; thus $fh = 0$
    and $(1-f)h = h$.

      Since $h \in C(X), \, \phi(h)= K_h +\tau(h)I
    \in \mathcal{K}^+$.  Note that $\tau(h) > 0$, since $h > 0$ on some neighborhood of
    $s_2$, and that the essential spectrum of $\phi(h)$ is
    $\{\tau(h)\}$. Now compute
        \begin{align*}\psi_1(h) &= Q_1^{-\frac{1}{2}}\left(\alpha \phi(fh) +
        \beta \phi((1-f)h)\right) Q_1^{-\frac{1}{2}}\\
        &= \beta Q_1^{-\frac{1}{2}} \phi(h) Q_1^{-\frac{1}{2}}\\
        & = \beta Q_1^{-\frac{1}{2}}K_hQ_1^{-\frac{1}{2}} + \beta
        \tau(h) Q_1^{-1}
        \end{align*}
    The first term in this sum is compact, while the second
    term can be written
        $$\beta \tau(h) Q_1^{-1} = \beta  \tau(h)[(\alpha -
        \beta)K_f + ((\alpha - \beta) \tau(f) + \beta)I]^{-1},$$
    where $\phi(f) = K_f + \tau(f)I$.
    Thus
        $$(q \circ \psi_1)(h) = \frac{\beta  \tau(h)}{(\alpha - \beta) \tau(f) +
        \beta}I + \mathcal{K}.$$
    Similar computations yield
        \begin{align*}
        \psi_2(h) &= (1 - \beta)Q_2^{-\frac{1}{2}} K_h
        Q_2^{-\frac{1}{2}} + (1 - \beta)  \tau(h) Q_2^{-1},
        \text{ and}\\
        (q \circ \psi_2)(h) &= \frac{(1 - \beta) \tau(h)}{(\beta - \alpha) \tau(f) + (1
        -\beta)}I + \mathcal{K}.
        \end{align*}
    So the essential spectra of $\psi_1(h)$ and $\psi_2(h)$ are
        $$\left\{\frac{\beta  \tau(h)}{(\alpha - \beta) \tau(f) +
        \beta}\right\} \text{ and } \left\{\frac{(1 - \beta) \tau(h)}{((\beta - \alpha) \tau(f) + (1
        -\beta)}\right\},$$
    respectively.  However, if these are equal, then
        $$\beta(\beta - \alpha) \tau(f) + \beta(1-\beta) = (1-\beta)(\alpha - \beta) \tau(f)+ \beta(1 - \beta),$$
    so that,
        $$\beta = \beta - 1,$$
    which is clearly impossible.  This shows that the essential spectra of $\psi_1(h)$ and
    $\psi_2(h)$ are distinct, so that
    $\psi_1(h)$ and $\psi_2(h)$ are not unitarily equivalent.  Thus
        $$\phi = t_1^* \psi_1t_1 + t_2^* \psi_2t_2$$
    expresses $\phi$ as a proper $C^*$-convex combination of
    positive unital maps $\psi_1$ and $\psi_2$ which are not both
    unitarily equivalent to $\phi$, demonstrating that $\phi$ is not
    $C^*$-extreme.  This proves the lemma.

\end{proof}

We can now prove the following:
\begin{theorem}\label{cmpct}
Let $\phi:C(X) \rightarrow \mathcal{K}^+$ be unital and positive.
Then $\phi$ is $C^*$-extreme if and only if $\phi$ is a
homomorphism.
\end{theorem}

\begin{proof}
If $\phi$ is multiplicative, then $\phi$ is $C^*$-extreme
\cite[Proposition 1.2]{FarenickMorenzC*ExPoGeStSpC*Al}. Conversely,
if $\phi$ is $C^*$-extreme, Lemma \ref{taumult} shows that the map
$\tau = q \circ \phi$ is multiplicative, so $\tau$ is a point evaluation $\tau(f) = f(s_0)$
for some point $s_0 \in X$.

Let $N$ be any neighborhood of $s_0$.  Then there exists a
continuous function $g_{_N}:X \rightarrow [0,1]$ with $g_{_N}(s_0) =
0$ and $g_{_N}|_{N^C} = 1$.

In this case $\tau(g_{_N})=0$, so
    $$\phi(g_{_N}) = K_{g_{_N}} \in \mathcal{K}.$$
Note that $\chi_{_{N^C}} \leq g_{_N}$, so that $\phi(\chi_{_{N^C}})
\leq \phi(g_{_N})$.
 Since $\mathcal{K}$ is hereditary, it follows that
 $\phi(\chi_{_{N^C}})$ is compact.  By Theorem~\ref{nec}, either
$\phi(\chi_{_{N^C}})$ is a projection or
$\sigma(\phi(\chi_{_{N^C}})) = [0,1]$. As a compact operator cannot
have the unit interval as its spectrum, $\phi(\chi_{_{N^C}})$ must
be a projection of finite rank.  Thus $\phi(\chi_{_N})$ is also a
projection.

Let $B$ be any Borel set of $X$ which does not contain $s_0$. Set
 $$\Lambda := \{K \subseteq B: K \text{ closed} \},$$
and partially order $\Lambda$ by inclusion. Then $\mu_\phi(K)$ is an
increasing net of projections.  Thus the SOT-$\displaystyle{\lim_K
\mu_\phi(K) =:Q}$ exists, and is a projection, namely the projection
onto $\displaystyle{\bigcup_{K \in \Lambda}\text{ran }\mu_\phi(K)}$.
Since the measures $\mu_{x,x}$ are regular for any choice of $x \in
\mathcal{H}$, we have

    $$\mu_{x,x}(B) = \sup_{K \in \Lambda} \mu_{x,x}(K)$$
or, equivalently,
    \begin{align*}
    \langle\mu_\phi(B)x,x \rangle &=
    \sup_{K \in \Lambda} \langle \mu_\phi(K) x,x \rangle\\
    &= \langle Qx,x \rangle.
    \end{align*}
As this holds for any $x \in \mathcal{H}$,
    $$Q = \mu_\phi(B).$$
If $B$ is a Borel set in $X$ which does contain $s_0$, then the
preceeding argument shows that $\mu_\phi(B^C)$ is a projection. Thus
$\mu_\phi(B)$ is also a projection.  Hence $\mu_{\phi}$ is a
projection valued measure, and $\phi$ is a homomorphism.

\end{proof}

\begin{remark}
When $\phi:C(X) \rightarrow \mathcal{K}^+$, as in
Theorem~\ref{cmpct}, we can obtain more information regarding the
support of $\mu_\phi$.  We have shown above that for any closed set
$K$ with $s_0 \not\in K$, $\mu_\phi(K)$ is a finite rank projection,
say of rank $n$. If $s_1, s_2$ are distinct points of $K \cap \text{
supp }~\mu_\phi$, let $N_1 \subseteq K$ be a neighborhood of $s_1$
which does not contain $s_2$. Then $K \setminus N_1$ is closed and
$s_0 \not\in K \setminus N_1$, so $\mu_\phi(K\setminus N_1)$ is a
projection of finite rank and
    $$0 < \text{rank } \mu_\phi(K \setminus N_1)< \text{rank } \mu_\phi (K) = n.$$
Since
    $$\mu_\phi(K) = \mu_\phi(K \setminus N_1) + \mu_\phi(N_1),$$
it follows that $\mu_\phi(N_1)$ is also a projection with
    $$0 < \text{rank }\mu_\phi(N_1) < n.$$
Clearly this process can be iterated at most $n$ times; we conclude
that any closed set $K \not\ni s_0$ contains at most finitely many
points of supp $\mu_\phi$.  Consequently, supp~$\mu_\phi$~$\setminus
\{s_0\}$ is a discrete set with at most one accumulation point at
$s_0$.
\end{remark}

If $\mathcal{H}$ is a separable Hilbert space, then it is clear from
the proof of Theorem~\ref{cmpct} and the preceding remark that the
support of $\mu_\phi$ must be at most countable with a single limit
point at $s_0$. In this case, $\phi$ must have the form
    $$\phi(f) = \sum_{s \in \text{supp} (\mu_\phi)}f(s)P_s,$$
where $P_s = \mu_\phi(\{s\})$ is a finite rank projection for each
$s \neq s_0$.  The rank of $P_{s_{_0}}$, on the other hand, may be
finite or infinite.  The following example, in which we consider the
case of a nonseparable Hilbert space, illustrates the structure of
unital positive maps $\phi:C(X) \rightarrow \mathcal{K}^+$.

\begin{example}\label{nonsep}
Let $\mathcal{H}$ be a nonseparable Hilbert space with dimension at
least as great as the cardinality of $\mathbb{R}$, and let $X =
\mathbb{R} \cup\{\omega\}$ be the one point compactification of
($\mathbb{R},d)$, the reals equipped with the discrete topology.
 Choose an orthonormal set $\{e_s\}_{s \in \mathbb{R}}$ in
$\mathcal{H}$ indexed by the reals, and write $P_s$ for the
projection onto the span of $e_s$.  Then, for any function $f \in
C(X)$, the set
    $$S(f) := \{s \in X:f(s)\neq f(\omega)\}$$
is at most countable, and
    $$\lim_{n \rightarrow \infty} f(s_n) =f(\omega),$$
where $\{s_n\}$ is any enumeration of $S(f)$.  Define a positive map
$\phi$ on $C(X)$ by
    $$\phi(f) = \sum_{s \in S(f)}[f(s) - f(\omega)]P_s + f(\omega)I.$$
Then for each $s \in \mathbb{R}$, the function $\delta_s =
\chi_{_{\{s\}}}$ is continuous and $\phi(\delta_s) = \mu_\phi(\{s\})
= P_s$. As in the proof of Theorem~\ref{cmpct}, if $G$ is any
neighborhood of $\omega$, then $G^C$ is a closed set not containing
$\omega$, and  $\phi(\chi_{_G})$ is a projection.  In this case the
descending net $\phi(\chi_{_G})$ of projections converges to the
projection $\phi(\chi_{\{\omega\}}) = 0$.  Thus $\mu_\phi$ is a
projection valued measure.

Note that we could define similar maps $\phi_1$ and $\phi_2$ by
    $$\phi_1(f) = \sum_{s \in S(f)}[f(s) - f(\omega)]P_{1/s} + f(\omega)I,$$
and
    $$\phi_2(f) = \sum_{s \in S(f)}[f(s) - f(\omega)]P_{\arctan s} + f(\omega)I.$$
For these two maps, we have $\phi_1(\chi_{_{\{\omega\}}})= P_0$,
while $\phi_2(\chi_{_{\{\omega\}}})$ is the projection onto the
closed span $\{\text{ran }P_s:s \in (-\infty, \pi /2] \cup [\pi /2,
\infty)\}$. Thus, the measure of $\{\omega\}$ may be a projection of
either finite or infinite rank.

\end{example}

\bibliographystyle{plain}
\bibliography{martha}

\end{document}